\documentclass{amsart}

\usepackage{amsmath}
\usepackage{amscd}
\usepackage{amssymb}

\input xy
\xyoption{all}

\newcommand{\bA}{{\mathbb A}}
\newcommand{\bC}{{\mathbb C}}

\newcommand{\bF}{{\mathbb F}}

\newcommand{\bP}{{\mathbb P}}
\newcommand{\bQ}{{\mathbb Q}}

\newcommand{\bZ}{{\mathbb Z}}

\newcommand{\cO}{{\mathcal O}}

\newcommand{\cX}{{\mathcal X}}

\DeclareMathOperator{\age}{age}
 
\DeclareMathOperator{\diag}{diag}

\DeclareMathOperator{\Spec}{Spec}

\newtheorem{theorem/definition}{Theorem/Definition}[section]

\theoremstyle{remark}

\theoremstyle{definition}

\newcommand{\be}{\begin{equation}}
\newcommand{\ee}{\end{equation}}
\newcommand{\bea}{\begin{eqnarray}}
\newcommand{\eea}{\end{eqnarray}}
\newcommand{\ben}{\begin{eqnarray*}}
\newcommand{\een}{\end{eqnarray*}}
\newcommand{\bet}{\begin{equation}
\begin{split}}
\newcommand{\eet}{\end{split}
\end{equation}}

\begin{document}

\title
{Arithmetic McKay Correspondence}
\author{Jian Zhou}
\address{Department of Mathematical Sciences\\Tsinghua University\\Beijing, 100084, China}
\email{jzhou@math.tsinghua.edu.cn}

\begin{abstract}
We propose an arithmetic McKay correspondence
which relates suitably defined zeta functions of some Deligne-Mumford stacks
to the zeta functions of their crepant resolutions.
Some examples are discussed.
\end{abstract}
\maketitle

\section{Introduction}

Let $\cX$ be an orbifold (smooth Deligne-Mumford stack)
with coarse moduli space $X=|\cX|$ and let $\pi: Y \to X$ be a crepant resolution.
As suggested by string theory on orbifolds \cite{Dix-Har-Vaf-Wit, Zas},
invariants of $Y$ should coincide with suitably defined orbifold invariants of $\cX$
 when $\cX$ satisfies the Gorenstein condition.
This is often referred to as the McKay correspondence in the mathematical literature.
See Reid \cite{Rei} for an exposition of some examples for classical invariants,
e.g. Euler numbers,
cohmomology groups, K-theory and derived categories, etc.
The correspondence between the Gromov-Witten invariants of $Y$ and the orbifold Gromov-Witten invariants
of $\cX$ is called the quantum McKay correspondence.
For example,
the Crepant Resolution Conjecture of Ruan \cite{Rua} states that the small orbifold quantum cohomology of $\cX$
is related to the small quantum cohomology of $Y$.
There have recently appeared much work and some generalizations of this conjecture,
see e.g. \cite{Bry-Gra, Zho2} and the references therein.

The McKay correspondence and the quantum McKay correspondence are usually considered
over the field of complex numbers.
It is very natural to consider their arithmetic versions by considering
good reductions  $\cX(\bF_q)$ over a finite field $\bF_q$.
Motivated by the Weil Conjecture and McKay correspondence,
we will define an orbifold zeta function $Z_{orb}(\cX, t)$ and
conjecture that
\be
Z_{orb}^\cX(t) = Z^Y(t),
\ee
where $Z^Y(t)$ is the zeta function of $Y$.
We will refer to this equality as the arithmetic McKay correspondence.
We will provide some examples in this paper.

During the preparation of this work,
we notice that Rose \cite{Ros} defines an $l$-adic cohomology version of the Chen-Ruan cohomology ring
and defines an orbifold cohomological zeta function using the orbifold arithmetic
Frobenius action.
He makes a similar conjecture using his zeta function.

\section{Arithmetic McKay Corresondence}

\subsection{Zeta functions of algebraic varieties}
Let $k = \bF_q$ be the finite field with $q$ elements,
and let $\bar{k}$ be its algebraic closure.
Suppose that $X$ is a scheme of finite type over $k$,
and let $\bar{X} := X \times_{\Spec k} \Spec \bar{k}$.
For every integer $r \geq 1$,
denote by $N(X(\bF_{q^r}))$ the number of $\bF_{q^r}$-rational points of $\bar{X}$.
Define the zeta function of $X$ by
\be
Z^X(t) = \exp \big( \sum_{r\geq 1} N(X(\bF_{q^r})) \frac{t^r}{r}\big).
\ee
By the Weil Conjecture \cite{Del},
for $X$ a smooth projective variety of dimension $n$ over $k = \bF_q$,
the zeta function $Z^X(t)$ has the following properties:
\begin{itemize}
\item[(a)] (Rationality) $Z^X(t)$ is a rational function of $t$:
\be
Z^X(t) = \frac{P_1(t) P_3(t) \cdot P_{2n-1}(t)}{P_0(t) P_2(t) \cdots P_{2n}(t)},
\ee
where $P_0(t) = 1 - t$, $P_{2n}(t) = 1 - q^nt$,
and for $1 \leq i \leq 2n-1$,
$P_i(t)$ is a polynomial in $\bZ[t]$,

\item[(b)] (Betti numbers)
If $X$ is obtained from the reduction of a projective variety $Y$ over $\bC$,
then one has for $1 \leq i \leq 2n-1$
\be
\deg P_i(t) = b_i(Y),
\ee
where $b_i(Y)$ is the $i$-th Betti number of $Y$.

\item[(c)] (Functional equation)
Let $\chi(X)$ be the self-intersection number of the diagonal $\Delta \subset X \times X$.
Then $Z^X(t)$ satisfies the following equation:
\be
Z^X(\frac{1}{q^nt}) = \pm q^{n \chi(X)/2} t^{\chi(X)} Z^X(t).
\ee

\item[(c)] (Analogue of the Riemann hypothesis)
The polynomials $P_i(t)$ can be written as:
\be
P_i(t) = \prod_{j=1}^{b_i} (1 - \alpha_{ij}t),
\ee
where  the $\alpha_{ij}$'s are algebraic integers with $|\alpha_{ij}| = q^{i/2}$.
\end{itemize}

Some of these properties also hold for singular varieties.

\subsection{Zeta functions of Gorentstein Deligne-Mumford stacks}

Let $\cX$ be a smooth Deligne-Mumford  stack with coarse moduli space $X$ over $\bF_q$.
Denote by $I\cX = \coprod_i I_i\cX$ its inertia stack,
where each $I_i\cX$ is a connected component.
Suppose that $q$ is coprime to the order of the automorphism group of each object in $\cX$.
Then one can define the {\em age} function $\age: I\cX \to \bQ$.
We say $\cX$ is Gorenstein if the age function takes integral values.

For a smooth Gorenstein Deligne-Mumford stack $\cX$ over $F_q$,
define
\be
N_{orb}(\cX(\bF_{q})) = \sum_i q^{\age(I_i(\cX))} N(|I_i\cX|(\bF_{q})),
\ee
where $|I_i\cX|$ is the coarse moduli space of $I_i(\cX)$.
Define the {\em orbifold zeta function} $Z_{orb}(\cX)$ by:
\be
Z_{orb}^\cX(t) = \exp(\sum_{r \geq 1} N_{orb}(\cX) \frac{t^r}{r}).
\ee

This zeta function contains information about the orbifold Betti numbers of $\cX(\bC)$
by the Lefschetz trace formula for Deligne-Mumford stacks \cite{Beh}.

\subsection{Arithmetic McKay correspondence}

Suppose that $\cX$ is a smooth Gorenstein Deligne-Mumford stack $\cX$ over $F_q$,
such that there is a crepant resolution $\pi: Y \to X = |\cX|$,
then we conjecture that:
\be
N_{orb}(\cX(\bF_q)) = N(Y(\bF_q)),
\ee
and so
\be
Z_{orb}^{\cX}(t) = Z^Y(t),
\ee
when $q$ is coprime to the order of the automorphism groups of objects in $\cX$.
We will refer to this equality as the {\em arithmetic McKay correspondence}.
By the Lefschetz trace formula for Deligne-Mumford stacks \cite{Beh},
we also conjecture that the smooth cohomology $H^*(\overline{\cX}_{sm}, \bQ_l)$ of $\overline{\cX}$
is isomorphic to the \'etale cohomology $H^*_{et}(|\cX|, \bQ_l)$ of the coarse moduli space $|\cX|$,
and the arithmetic Frobenius actions on these spaces can be identified.
We hope our conjecture sheds some light on the results of  \cite{Bat, Den-Loe, Lup-Pod, Yas}.

\section{Some Examples}

\subsection{The case of Kleinian singularities}
\label{sec:Kleinian}

It is well-known that the binary polyhedral groups $G \subset SL(2)$ are finite algebraic groups.
For their realizations and invariants see \cite{Spr}.
We will now verify the arithmetic McKay correspondence for the stacks $[\bA^2/G]$.
The coarse moduli space $|\bA^2/G|$ has a Kleinian singularity:
\begin{align}
\text{Cyclic group} && xy=z^n, \\
\text{Binary dihedral group} && x^2 + y^2z + z^{n+1} = 0, \\
\text{Binary tetrahedral group} && x^2+y^3+z^4 = 0, \\
\text{Binary octahedral group} && x^2 +y^3z + z^3 = 0, \\
\text{Binary icosahedral group} && x^2 + y^3 + z^5 = 0.
\end{align}
Now $I[\bA^2/G]$ has $k$ twisted components, each has only one object and each is of age $1$,
where $k$ is the number of nontrivial conjugacy classes.
Therefore,
\be \label{eqn:Orb}
N_{orb}([\bA^2/G](\bF_q)) = N(|\bA^2/G|(\bF_q)) + k q.
\ee
On the other hand,
the minimal resolution $\widehat{|\bA^2/G|}$ is a crepant resolution,
and the exceptional divisor is a nodal curve with $k$ components which are isomorphic to $\bP^1$;
furthermore, there are exactly $k-1$ nodal points.
Therefore,
\ben
N(\widehat{|\bA^2/G|}(\bF_q))
& = & N(|\bA^2/G|(\bF_q)) - 1 + k N(\bP^1(\bF_q)) - (k-1) \\
& = &  N(|\bA^2/G|(\bF_q)) -1 + k (q + 1) - (k-1) \\
& = & N(|\bA^2/G|(\bF_q)) + k q.
\een
Hence by comparing with (\ref{eqn:Orb}),
we have
\be
N(\widehat{|\bA^2/G|}(\bF_q))
= N_{orb}([\bA^2/G](\bF_q)).
\ee

We conjecture that
\be
N(\widehat{|\bA^2/G|}(\bF_q))
= N_{orb}([\bA^2/G](\bF_q)) = q^2 + k q,
\ee
i.e., $N(|\bA^2/G|(\bF_q)) = q^2$,
and so
\be
Z_{orb}^{[\bA^2/G]}(t) = Z^{\widehat{|\bA^2/G|}}(t) = \frac{1}{(1-q^2t)(1-qt)^k}.
\ee
We can check the following cases.
For the cyclic group case,
the equation
$$xy=z^n$$
clearly has $q^2$ solutions in $\bF^2_q$.
Indeed,
for $z=0$,
we have $2q-1$ solutions: $x=0$, $y \in \bF_q$ arbitrary or $x \in \bF_q^*$, $y=0$;
for $z\in \bF_q^*$,
we have $(q-1)^2$ solutions: $y,z \in \bF_q^*$ arbitrary and
$$x = \frac{1}{y} z^n.$$
For the binary tetrahedral group case,
when $3 \not| (q-1)$,
there is a unique solution to each equation $y^3 = a$ ($a\in \bF_q$)
because $\bF_q^*$ is a cyclic group of order $q-1$.
Therefore,
we can take $x, z \in \bF_q$ arbitrary,
and take $y$ to be the unique solution to $y^3 = - x^2 - z^4$.
Similarly,
in the binary icosahedral case,
when $3 \not| (q-1)$ or $5 \not| (q-1)$,
the equation $x^2 + y^3+z^5 = 0$ has $q^2$ solutions.

\subsection{The case of Georenstein singularities in dimension $3$}

Let $G \subset SL(3)$ be a finite group.
They have been classified over $\bC$ by Miller {\em et al.} \cite{Mil-Bli-Dic}.
For the crepant resolution $\widehat{|\bA^3/G|}$ of $|\bA^3/G|$ in the complex case see Roan \cite{Roa}
and the references therein.
We expect that it is possible check the arithmetic McKay correspondence using explicit constructions
of the crepant resolutions.
For example,
let $\mu_3$ act on $\bA^3$ by multiplication by $\xi_3$, a primitive root of unity of order $3$.
There are two twisted sectors,
one with age $1$ and one with age $2$.
Therefore,
\be
N_{orb}([\bA^3/\mu_3](\bF_q)) = N(|\bA^3/\bZ_3|(\bF_q)) + q + q^2.
\ee
The exceptional divisor of the crepant resolution of $|\bA^3/\mu_3|$ is $\bP^2$,
therefore,
\be
N(\widehat{\bA^3/\mu_3}(\bF_q))
= N(|\bA^3/\mu_3|(\bF_q)) - 1 + N(\bP^2(\bF_q))
= N(|\bA^3/\mu_3|(\bF_q)) + q + q^2.
\ee
Hence we get:
\be
N_{orb}([\bA^3/\mu_3](\bF_q)) = N(\widehat{\bA^3/\mu_3}(\bF_q)).
\ee
For another example,
let a generator of $\mu_5$ act on $\bA^3$ by the diagonal matrix
$\diag( \xi_5, \xi_5^2, \xi_5^2)$.
Then $[\bA^3/\mu_5]$ has four twisted sectors, with ages $1$, $1$, $2$ and $2$, respectively.
Therefore,
\be
N_{orb}([\bA^3/\mu_5](\bF_q)) = N(|\bA^3/\bZ_3|(\bF_q)) + 2q + 2q^2.
\ee
The exceptional divisor of the crepant resolution of $|\bA^3/\mu_5|$ is a copy $\bP^2$ glued to a copy
of the Hirzebruch surface $\bP(\cO_{\bP^1} \oplus \cO_{\bP^1}(3))$ along with a copy of $\bP^1$,
therefore,
\ben
&& N(\widehat{\bA^3/\mu_5}(\bF_q)) \\
& = & N(|\bA^3/\mu_5|(\bF_q)) - 1 + N(\bP^2(\bF_q)) + N(\bP(\cO_{\bP^1} \oplus \cO_{\bP^1}(3))(\bF_q))
- N(\bP^1(\bF_q)) \\
& = & N(|\bA^3/\mu_5|(\bF_q)) - 1 + (1 + q + q^2) + (1+q)^2 - (1+q) \\
& = & N(|\bA^3/\mu_5|(\bF_q)) + 2 q + 2 q^2.
\een
Here we have used the following fact:
For a fiber bundle $X \to B$ with fiber $F$ over $\bF_q$,
we have
\be
N(X(\bF_q)) = N(B(\bF_q)) \cdot N(F(\bF_q)).
\ee
Hence we get:
\be
N_{orb}([\bA^3/\mu_5](\bF_q)) = N(\widehat{\bA^3/\mu_5}(\bF_q)).
\ee

Similar to the two dimensional case,
we conjecture the following equality
\be
N(|\bA^3/G|(\bF_q)) = q^3
\ee
holds for a finite group $G \subset SL(3)$.

\subsection{Symmetric products and Hilbert schemes of algebraic surfaces}
Let $X$ be a smooth projective surface defined over $\bF_q$.
Denote by $X^{(n)} = |X^n/S_n| $ the $n$-th symmetric product of $S$ and
by $S^{[n]}$ the Hilbert schemes of $0$-dimensional subscheme of length $n$.
According to Forgarty \cite{For}, the Hilbert-Chow morphism $\pi: X^{[n]} \to X^{(n)}$ is a crepant resolution.
The components of $I(X^n/S_n)$ are parameterized by partitions $\mu:=1^{m_1} \cdots n^{m_n}$ of $n$,
where $m_1, \dots, m_n$ are nonnegative integers such that
$$\sum_{i=1}^n m_i i = n.$$
Denote by $X^{(n)}_\mu$ this component.
It is isomorphic to $\prod_{i=1}^n X^{(m_i)}$,
with age $a_{\mu}:=\sum_{i=1}^n m_i(i-1)$ (see e.g. \cite{Zho1}).
According \cite[p. 196]{Got},
\ben
\sum_{n=0}^{\infty} N(X^{(n)}(\bF_q)) t^n
= Z^X(t)
= \exp \sum_{r=1}^{\infty} N(X(\bF_{q^r})) \frac{t^r}{r}.
\een
Therefore,
we have
\ben
&& \sum_{n=0}^\infty Q^n N_{orb}(X^{(n)}(\bF_q))
= \sum_{n \geq 0} Q^n \sum_{\mu \vdash n} q^{a_\mu} N(X^{(n)}_{\mu}(\bF_q)) \\
& = & \sum_{n=0}^\infty Q^n \sum_{\sum_{i=1}^n m_i i =n}
q^{\sum_{i=1}^n m_i(i-1)} \prod_{i=1}^n N(X^{(m_i)}(\bF_q)) \\
& = & \prod_{i=1}^\infty \sum_{m_1 = 0}^\infty
(q^{i-1}Q^i)^{m_i} N(X^{(m_i)}(\bF_q))
= \prod_{i=1}^\infty Z^X(q^{i-1}Q).
\een
The key technical result in \cite{Got} is (Lemma 2.9 on p. 203):
\be
\sum_{n=0}^\infty Q^n N(X^{[n]}(\bF_q)) = \prod_{i=1}^\infty Z^X(q^{i-1}Q).
\ee
This is equivalent to the arithmetic McKay correspondence in this case:
$$N_{orb}(X^{(n)}(\bF_q)) = N(X^{[n]}(\bF_q)).$$

Suppose that $Y$ is a projective surface acted on by a finite group $G$ such that the
singularities of $|Y/G|$ are Kleinian singularities.
Let $\pi: X \to Y$ be its crepant resolution.
Now the wreath product group $G_n$ acts on $Y^n$ and we have a crepant resolution
$X^{[n]} \to |Y^n/G_n|$.
The inertia stack of $[Y^n/G_n]$ is described in \cite[\S 2.1]{Wan-Zho},
and the ages are computed in \cite[\S 3.2]{Wan-Zho}.
It is straightforward to combine G\"ottsche's result and the results for Kleinian singularities
in \S \ref{sec:Kleinian} to verify the arithmetic McKay correspondence in this case:
$$N_{orb}(|Y^n/G_n|(\bF_q)) = N(X^{[n]}(\bF_q)).$$
The details are left to the interested reader.

{\em Acknowledgements}.
This research is partially supported by two NSFC grants (10425101 and 10631050)
and a 973 project grant NKBRPC (2006cB805905).

\end{document}